\documentclass[12pt]{article}%
\usepackage{amssymb}
\usepackage{amsmath}
\usepackage{amsfonts}
\usepackage{cite}
\usepackage{graphicx}
\usepackage{subfig}%

\begin{document}
\author{I.P. Smirnov\\Institute of Applied Physics RAS, \\46 Ul'yanova Street, Nizhny Novgorod, Russia}
\title{The reduction of the problem of maximization of the fraction of two functionals }
\date{}
\maketitle

\begin{abstract}
We propose an algorithm for reduction of the problem of maximization of
fraction of two functionals to the equivalent procedure including maximization
of difference between the functionals and the solution of an equation of
scalar unknown. For illustration of the algorithm we solve some problems of
the described type.
\end{abstract}


\section{Introduction}

\textbf{Theorem 1}. \textit{Let }$W\left(  x\right)  >0,$\textit{ }%
$W_{0}\left(  x\right)  $\textit{ be continuous functionals on a compact set
}$K$\textit{ of a metric space, }$\beta$\textit{ be a number,}%
\begin{align*}
J\left(  x\right)   &  \equiv\frac{W_{0}\left(  x\right)  }{W\left(  x\right)
},\\
J_{\beta}\left(  x\right)   &  \equiv W_{0}\left(  x\right)  -\beta W\left(
x\right)  .
\end{align*}
\textit{Consider the following extremal problems}%
\begin{equation}
J\rightarrow\max\limits_{K}, \label{Jmax}%
\end{equation}
\textit{and}%
\begin{equation}
J_{\beta}\rightarrow\max\limits_{K}. \label{Jbmax}%
\end{equation}
\textit{Let }$x_{\max}$ \textit{be a solution of the problem (\ref{Jmax}),
}$x_{\beta}$\textit{ a solution of (\ref{Jbmax}),}%
\[
\beta_{\max}\equiv J\left(  x_{\max}\right)  ,\ \ \ j\left(  \beta\right)
\equiv J_{\beta}\left(  x_{\beta}\right)  .
\]
\textit{Then}%
\begin{equation}
\left\{
\begin{array}
[c]{cc}%
j\left(  \beta\right)  <0, & \beta>\beta_{\max},\\
j\left(  \beta\right)  >0, & \beta<\beta_{\max},\\
j\left(  \beta\right)  =0, & \beta=\beta_{\max},
\end{array}
\right.  \label{betamax}%
\end{equation}
\textit{so }$\beta_{\max}$\textit{ is the only solution of the equation}%
\begin{equation}
j\left(  \beta\right)  =0. \label{jbeta}%
\end{equation}
\textit{ Functionals }$J\left(  x\right)  $\textit{ and }$J_{\beta_{\max}%
}\left(  x\right)  $\textit{ take their maxima at the same point }$x_{\max
}=x_{\beta_{\max}}$\textit{. The function }$\beta\rightarrow$\textit{
}$x_{\beta}$\textit{ is continuous in every compact segment }$\left[
\beta_{1},\beta_{2}\right]  $\textit{.}

\textbf{Prove.} Under the hypotheses of the theorem, the solutions of the
extrema problems $x_{\max},$\ $x_{\beta}$ exist (may be not unique). If
$\beta>\beta_{\max}$, then $J_{\beta}\left(  x\right)  =W\left(  x\right)
\left(  J\left(  x\right)  -\beta\right)  <0$ for all $x$; so $j\left(
\beta\right)  =J_{\beta}\left(  x_{\beta}\right)  <0$. If $\beta<\beta_{\max}%
$, then $J_{\beta}\left(  x_{\max}\right)  =W\left(  x_{\max}\right)  \left(
J\left(  x_{\max}\right)  -\beta\right)  =W\left(  x_{\max}\right)  \left(
\beta_{\max}-\beta\right)  >0,$ so $j\left(  \beta\right)  =J_{\beta}\left(
x_{\beta}\right)  \geq J_{\beta}\left(  x_{\max}\right)  >0$. For $\beta
=\beta_{\max}$ we have%
\[
J_{\beta_{\max}}\left(  x\right)  =W\left(  x\right)  \left(  J\left(
x\right)  -\beta_{\max}\right)  \leq0,
\]
for all $x$ and%
\[
J_{\beta_{\max}}\left(  x_{\max}\right)  =W\left(  x_{\max}\right)  \left(
J\left(  x_{\max}\right)  -\beta_{\max}\right)  =0,
\]
so $j\left(  \beta_{\max}\right)  =\max_{x}J_{\beta_{\max}}\left(  x\right)
=0$ and we can take $x_{\beta_{\max}}=x_{\max}$. The last statement of the
theorem is the implication of the continuity of the functional $J_{\beta}$ on
the compact $\left[  \beta_{1},\beta_{2}\right]  \times K$ $\blacktriangle$

\textbf{Theorem 2}. \textit{Let }$W\left(  x\right)  \neq0$ and\textit{
}$W_{0}\left(  x\right)  $\textit{ be functionals on a set }$K$\textit{ of a
metric space}%
\begin{align*}
J\left(  x\right)   &  \equiv\frac{W_{0}\left(  x\right)  }{W\left(  x\right)
},\\
J_{\beta}\left(  x\right)   &  \equiv W\left(  x\right)  \left(  W_{0}\left(
x\right)  -\beta W\left(  x\right)  \right)  .
\end{align*}
\textit{Suppose that for every }$\beta$ \textit{the extremal problems}%
\[
J\rightarrow\max\limits_{K},\ \ J_{\beta}\rightarrow\max\limits_{K}%
\]
\textit{have their solutions }$x_{\max},x_{\beta}$\textit{. Let }%
\[
\beta_{\max}\equiv J\left(  x_{\max}\right)  ,\ \ \ j\left(  \beta\right)
\equiv J_{\beta}\left(  x_{\beta}\right)  .
\]
\textit{Then}%
\[
\left\{
\begin{array}
[c]{cc}%
j\left(  \beta\right)  <0, & \beta>\beta_{\max},\\
j\left(  \beta\right)  >0, & \beta<\beta_{\max},\\
j\left(  \beta\right)  =0, & \beta=\beta_{\max},
\end{array}
\right.
\]
\textit{so }$\beta_{\max}$\textit{ is the only solution of the equation}%
\[
j\left(  \beta\right)  =0.
\]
\textit{ Functionals }$J\left(  x\right)  $\textit{ and }$J_{\beta_{\max}%
}\left(  x\right)  $\textit{ take their maxima at the same point }$x_{\max
}=x_{\beta_{\max}}$\textit{. }

\textbf{Prove.} To prove the theorem we multiply the numerator and denominator
of the fraction $J$ by $W\left(  x\right)  $ and repeat the prove of the
previous theorem $\blacktriangle$

The corollary of the theorems is the following procedure of solution of the
problem (\ref{Jmax}). On the first step we solve the problem (\ref{Jbmax}) for
any arbitrary $\beta$; then we calculate the function $j\left(  \beta\right)
$ and solve the scalar equation (\ref{jbeta}).

By calculation of $j\left(  \beta\right)  $ for any given $\beta$ we can see
the direction where root of the equation is situated. Then any iteration
scheme can be applied to find the root with the necessary accuracy.

\section{Examples}

Let us consider some problems where the above algorithm can be applied.

\begin{description}
\item[Problem 1.] Let $x\in$ $\left[  x_{1},x_{2}\right]  $, $W\left(
x\right)  =ax+b$, $W_{0}\left(  x\right)  =a_{0}x+b_{0}$, $ax_{1}+b>0$,
$ax_{2}+b>0$, $a>0$.

In this case%
\begin{align*}
J\left(  x\right)   &  \equiv\frac{W_{0}\left(  x\right)  }{W\left(  x\right)
}=\frac{a_{0}x+b_{0}}{ax+b},\ \ \ \ \\
J_{\beta}\left(  x\right)   &  \equiv W_{0}\left(  x\right)  -\beta W\left(
x\right)  =\left(  a_{0}-\beta a\right)  x+b_{0}-\beta b.
\end{align*}

It is clear that the solution of the problem (\ref{Jbmax}) is%
\[
x_{\beta}=\left\{
\begin{array}
[c]{cc}%
x_{1}, & \beta\geq a_{0}/a,\\
x_{2}, & \beta<a_{0}/a,
\end{array}
\right.
\]
so%
\[
j\left(  \beta\right)  =J_{\beta}\left(  x_{\beta}\right)  =\left\{
\begin{array}
[c]{cc}%
a_{0}x_{1}+b_{0}-\beta\left(  ax_{1}+b\right)  , & \beta\geq a_{0}/a,\\
a_{0}x_{2}+b_{0}-\beta\left(  ax_{2}+b\right)  , & \beta<a_{0}/a.
\end{array}
\right.
\]

Note that%
\[
j\left(  \frac{a_{0}}{a}\right)  =b_{0}-\frac{a_{0}b}{a}=\frac{ab_{0}-ba_{0}%
}{a}.
\]
The solution of the equation $j\left(  \beta\right)  =0$ is
\[
\beta_{\max}=J\left(  x_{\max}\right)  =\left\{
\begin{array}
[c]{cc}%
J\left(  x_{1}\right)  , & j\left(  \frac{a_{0}}{a}\right)
>0\Longleftrightarrow ab_{0}-ba_{0}>0,\\
J\left(  x_{2}\right)  , & j\left(  \frac{a_{0}}{a}\right)  \leq
0\Longleftrightarrow ab_{0}-ba_{0}\leq0.
\end{array}
\right.
\]
Hence, the solution of the problem (\ref{Jmax})%
\[
x_{\max}=\left\{
\begin{array}
[c]{cc}%
x_{1}, & ab_{0}-ba_{0}>0,\\
x_{2}, & ab_{0}-ba_{0}\leq0.
\end{array}
\right.
\]
We can also obtain the same result using the derivation
\[
J^{\prime}\left(  x\right)  =\frac{a_{0}b-ab_{0}}{\left(  ax+b\right)  ^{2}}.
\]

\item[Problem 2.] Let $x\in\left[  x_{1},x_{2}\right]  $, $W\left(  x\right)
=ax^{2}+bx+c$, $W_{0}\left(  x\right)  =a_{0}x^{2}+b_{0}x+c_{0}$, $ax_{1}%
^{2}+bx_{1}+c>0$, $ax_{2}^{2}+bx_{2}+c>0$, $a>0$.

We have now%
\begin{align*}
J\left(  x\right)   &  \equiv\frac{W_{0}\left(  x\right)  }{W\left(  x\right)
}=\frac{a_{0}x^{2}+b_{0}x+c_{0}}{ax^{2}+bx+c},\ \\
J_{\beta}\left(  x\right)   &  \equiv W_{0}\left(  x\right)  -\beta W\left(
x\right)  =\left(  a_{0}-\beta a\right)  x^{2}+\left(  b_{0}-\beta b\right)
x+c_{0}-\beta c.
\end{align*}

The solution of the problem $J_{\beta}\left(  x\right)  \rightarrow\max$ is
one of the following three values $x_{1},x_{2}$ and
\[
x_{3}=\frac{1}{2}\left(  b_{0}-\beta b\right)  /\left(  a_{0}-\beta a\right)
\]
(in the case when the last belongs to $\left[  x_{1},x_{2}\right]  $;
otherwise the solution is one of two values $x_{1},x_{2}$). Hence,%
\[
j\left(  \beta\right)  =\left\{
\begin{array}
[c]{ll}%
\max\left\{  J_{\beta}\left(  x_{1}\right)  ,J_{\beta}\left(  x_{2}\right)
,J_{\beta}\left(  x_{3}\right)  \right\}  , & x_{3}\in\left[  x_{1}%
,x_{2}\right]  ,\\
\max\left\{  J_{\beta}\left(  x_{1}\right)  ,J_{\beta}\left(  x_{2}\right)
\right\}  , & x_{3}\notin\left[  x_{1},x_{2}\right]  .
\end{array}
\right.
\]

\item[Problem 3.] Let $x\in\left[  x_{1},x_{2}\right]  $, $f_{0}\left(
x\right)  >0$, $f\left(  x\right)  >1$,%
\[
J\left(  x\right)  =\frac{\ln f_{0}}{\ln f}\rightarrow\max.
\]
The problem%
\[
\ln f_{0}-\beta\ln f\rightarrow\max
\]
is equivalent to%
\[
J_{\beta}\left(  x\right)  =\frac{f_{0}}{f^{\beta}}\rightarrow\max.
\]
For this problem we construct the auxiliary problem%
\[
J_{\gamma,\beta}=f_{0}-\gamma f^{\beta}\rightarrow\max.
\]
Let $x_{\gamma,\beta}$ be a solution of the last problem, $\gamma
=\gamma\left(  \beta\right)  $ a solution of the equation%
\[
J_{\gamma,\beta}\left(  x_{\gamma,\beta}\right)  =0
\]
for a fixed $\beta$. Next let $\beta_{\max}$ be a solution of the equation%
\[
J_{\beta}\left(  x_{\gamma\left(  \beta\right)  ,\beta}\right)  =1.
\]
Then $x_{\gamma\left(  \beta_{\max}\right)  ,\beta_{\max}}$ is a solution of
the initial problem.

\item[Problem 4.] Let $x\in K$, where $K=\left\{  x:\left\Vert x\right\Vert
\leq r\right\}  $ is a solid sphere of a Hilbert space $H$, $W\left(
x\right)  =\left\langle w,x\right\rangle +h$, $h>r\left\Vert w\right\Vert $,
$W_{0}=\left\langle w_{0},x\right\rangle +h_{0}$,%
\begin{align*}
J\left(  x\right)   &  \equiv\frac{W_{0}\left(  x\right)  }{W\left(  x\right)
}=\frac{\left\langle w_{0},x\right\rangle +h_{0}}{\left\langle
w,x\right\rangle +h},\ \ \ \ \\
J_{\beta}\left(  x\right)   &  \equiv W_{0}\left(  x\right)  -\beta W\left(
x\right)  =\left\langle w_{\beta},x\right\rangle +h_{\beta},\\
w_{\beta}  &  \equiv w_{0}-\beta w,\ h_{\beta}=h_{0}-\beta h.
\end{align*}

It is clear that the solution of the problem (\ref{Jbmax}) has now the form%
\begin{equation}
x_{\beta}=r\frac{w_{\beta}}{\left\Vert w_{\beta}\right\Vert },\label{xbeta}%
\end{equation}
therefore,%
\begin{equation}
j\left(  \beta\right)  =r\left\Vert w_{0}-\beta w\right\Vert +h_{0}-\beta
h.\label{jbeta1}%
\end{equation}
So, the problem of maximization of $J$ has been transformed to the solution of
the nonlinear equation (in unknown value $\beta)$%
\begin{equation}
r\left\Vert w_{0}-\beta w\right\Vert +h_{0}-\beta h=0.\label{urav}%
\end{equation}

Let us show that the curves $y=r\left\Vert w_{0}-\beta w\right\Vert $ and
$y=\beta h-h_{0}$ intersect. We have%
\begin{gather*}
\lim\limits_{\beta\rightarrow\infty}\frac{r\left\Vert w_{0}-\beta w\right\Vert
}{\beta}=r\left\Vert w\right\Vert \operatorname{sgn}\beta,\\
\lim\limits_{\beta\rightarrow\infty}\left[  r\left\Vert w_{0}-\beta
w\right\Vert -r\beta\operatorname{sgn}\beta\left\Vert w\right\Vert \right]
=\\
=\lim\limits_{\beta\rightarrow\infty}r\frac{\left\Vert w_{0}-\beta
w\right\Vert ^{2}-\left\Vert \beta w\right\Vert ^{2}}{\left\Vert w_{0}-\beta
w\right\Vert +\left\Vert \beta w\right\Vert }=\\
=r\lim\limits_{\beta\rightarrow\infty}\frac{\left\langle w_{0}-\beta
w,w_{0}-\beta w\right\rangle -\left\langle \beta w,\beta w\right\rangle
}{\left\Vert w_{0}-\beta w\right\Vert +\left\Vert \beta w\right\Vert }=\\
=r\lim\limits_{\beta\rightarrow\infty}\frac{\left\Vert w_{0}\right\Vert
^{2}-\beta\left\langle w_{0},w\right\rangle -\beta\left\langle w,w_{0}%
\right\rangle }{\left\Vert w_{0}-\beta w\right\Vert +\left\Vert \beta
w\right\Vert }=\\
=-r\operatorname{sgn}\beta\operatorname{Re}\left\langle w_{0},\tilde
{w}\right\rangle ,\ \tilde{w}\equiv w/\left\Vert w\right\Vert .
\end{gather*}
So,
\begin{equation}
r\left\Vert w_{0}-\beta w\right\Vert =r\left\Vert w\right\Vert \left\vert
\beta\right\vert -r\operatorname{sgn}\beta\operatorname{Re}\left\langle
w_{0},\tilde{w}\right\rangle +o\left(  1\right)  ,\ \beta\rightarrow
\infty.\label{asimp}%
\end{equation}
Under the condition $h_{0}+r\left\Vert w_{0}\right\Vert >0$ the curves
intersect for positive $\beta$. Under the condition $h_{0}+r\left\Vert
w_{0}\right\Vert \leq0$ they intersect for $\beta\leq0$. So, $\beta_{\max}>0$
in the first case and $\beta_{\max}\leq0$ in the second. The cause of it is
that the functional $J$ takes some positive values in the first case and only
non positive in the second.

Substiuting the asymptotes (\ref{asimp}) to the equation (\ref{urav}), we get
a priory valuations of maximal value of $J:$
\begin{equation}
\beta_{\max}=\max J\simeq\left\{
\begin{array}
[c]{cc}%
\frac{h_{0}-r\operatorname{Re}\left\langle w_{0},\tilde{w}\right\rangle
}{h-r\operatorname{Re}\left\langle w,\tilde{w}\right\rangle }, &
h_{0}+r\left\Vert w_{0}\right\Vert >0,\\
\frac{h_{0}+r\operatorname{Re}\left\langle w_{0},\tilde{w}\right\rangle
}{h+r\operatorname{Re}\left\langle w,\tilde{w}\right\rangle }, &
h_{0}+r\left\Vert w_{0}\right\Vert \leq0.
\end{array}
\right.  \label{abeta}%
\end{equation}
Note that this valuations are correct only for big values of $\left\vert
\beta_{\max}\right\vert $.

\begin{description}
\item[Example 1.] Let $H=\mathbf{R}^{10}$, $w_{0}=\left(
1,1,1,1,1,0,0,0,0,10\right)  $,

$w=\left(  1,0,0,0,0,1,1,1,1,1\right)  $, $r=1$, $h_{0}=15$, $h=2.7$.

In this case $h_{0}+r\left\Vert w_{0}\right\Vert \simeq25.25>0$. The solution
of the problem is presented on fig. \ref{fig:optvect1}. The value $\max
J\simeq43.61$. Graphs of the functions $\beta\rightarrow j\left(
\beta\right)  $ and $\beta\rightarrow J\left(  x_{\beta}\right)  $ in $\left[
0,\beta_{\max}\right]  $ are presented on fig. \ref{fig:jb1}. The process of
asymptotic estimations of the $\max J$ is illustrated by fig. \ref{fig:y1},
the formula (\ref{abeta}) for $\max J$ gives
\[
\max J\simeq\frac{h_{0}-r\operatorname{Re}\left\langle w_{0},\tilde
{w}\right\rangle }{h-r\operatorname{Re}\left\langle w,\tilde{w}\right\rangle
}\simeq41.95.
\]

\item[Example 2.] Let $H=\mathbf{R}^{10}$, $w_{0}=\left(
1,1,1,1,1,0,0,0,0,10\right)  $,

$w=\left(  1,0,0,0,0,1,1,1,1,1\right)  $, $r=1$, $h_{0}=-15$, $h=2.7$.

In this case $h_{0}+r\left\Vert w_{0}\right\Vert \simeq-4.75<0$. The solution
of the problem is presented on fig. \ref{fig:optvect2}. The value $\max
J\simeq-1.18$. Graphs of the functions $\beta\rightarrow j\left(
\beta\right)  $ and $\beta\rightarrow J\left(  x_{\beta}\right)  $ in $\left[
0,\beta_{\max}\right]  $ are presented on fig. \ref{fig:jb2}. The process of
asymptotic estimations of the $\max J$ is illustrated by fig. \ref{fig:y2},
the formula (\ref{abeta}) for $\max J$ gives%
\[
\max J\simeq\frac{h_{0}+r\operatorname{Re}\left\langle w_{0},\tilde
{w}\right\rangle }{h+r\operatorname{Re}\left\langle w,\tilde{w}\right\rangle
}\simeq-2.04.
\]

\begin{figure}[ptb]
\centering\subfloat[Example 1]{\label{fig:optvect1}\includegraphics[width=0.5\textwidth]{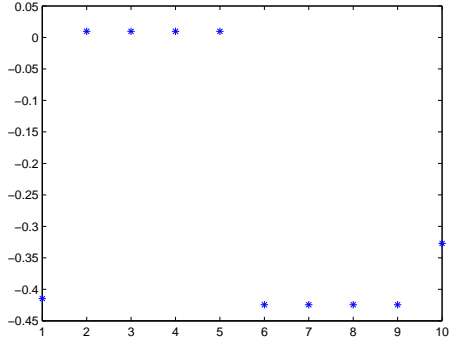}}
\subfloat[Example 2]{\label{fig:optvect2}\includegraphics[width=0.5\textwidth]{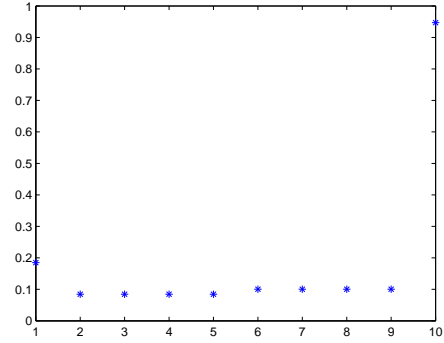}}\caption{Optimal
vectors $x_{\max}=x_{\beta_{\max}}$}%
\end{figure}\begin{figure}[ptb]
\centering\subfloat[Example 1]{\label{fig:jb1}\includegraphics[width=0.5\textwidth]{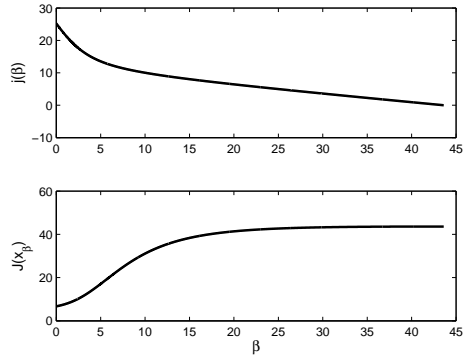}}
\subfloat[Example 2]{\label{fig:jb2}\includegraphics[width=0.5\textwidth]{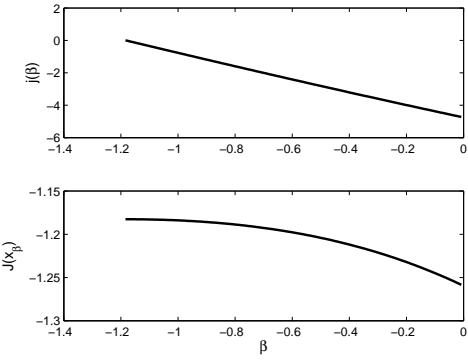}}\caption{Graphs
of functions $\beta\rightarrow j\left(  \beta\right)  $ (upper panels) and
functions $\beta\rightarrow J\left(  x_{\beta}\right)  $ (lower panels).}%
\end{figure}\begin{figure}[ptb]
\centering\subfloat[Example 1]{\label{fig:y1}\includegraphics[width=0.5\textwidth]{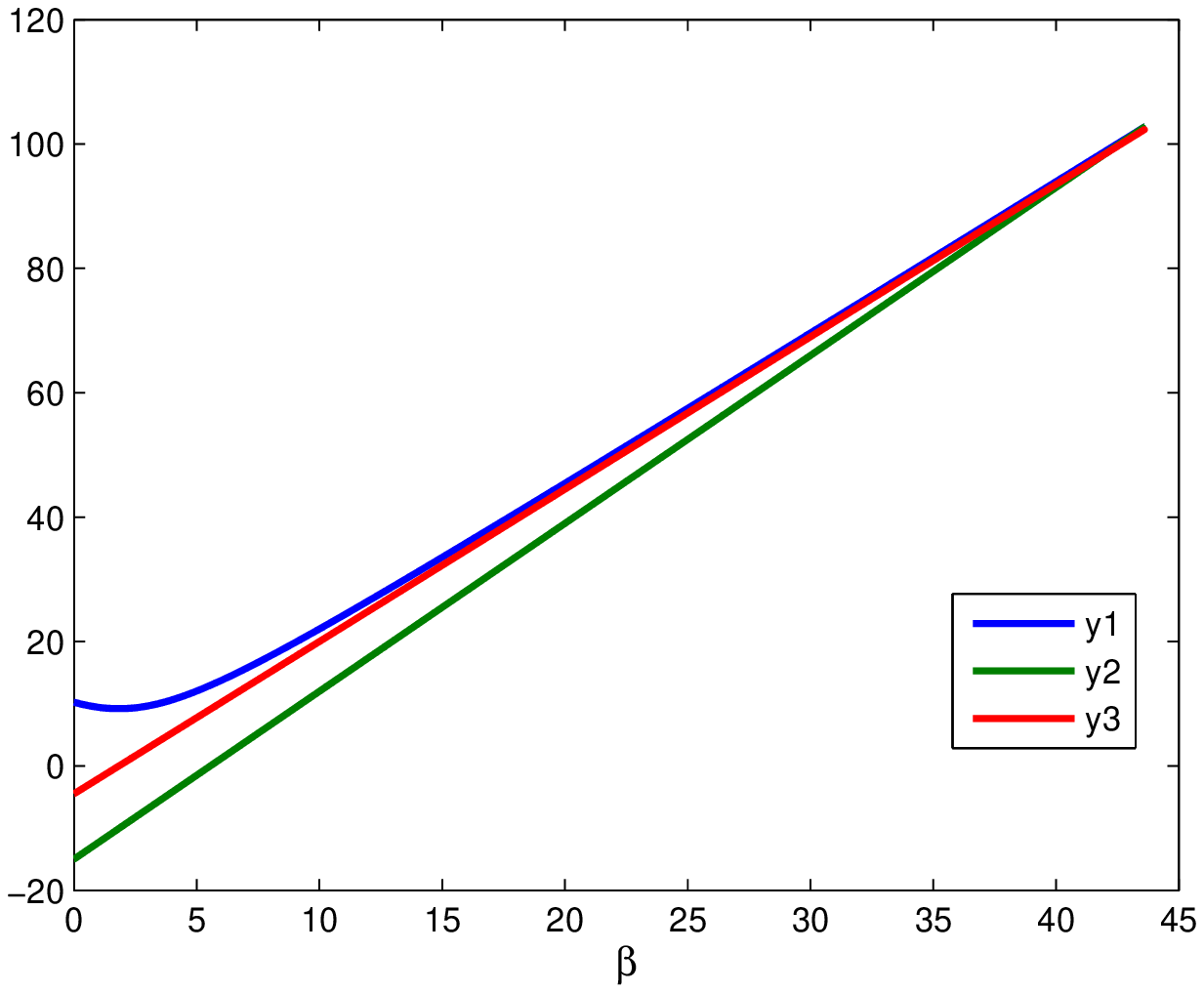}}
\subfloat[Example 2]{\label{fig:y2}\includegraphics[width=0.5\textwidth]{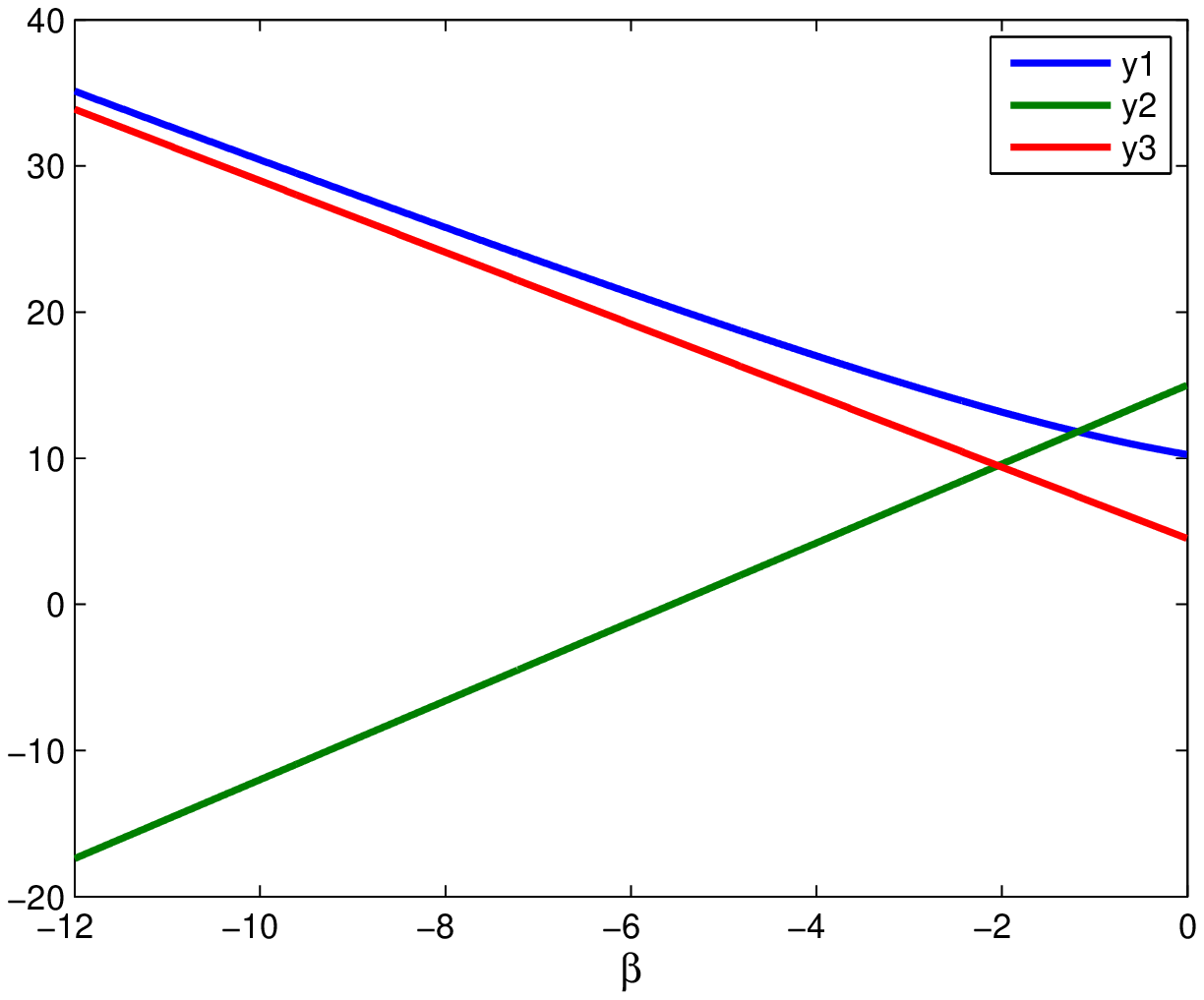}}\caption{Graphs
of the functions $y_{1}\left(  \beta\right)  =r\left\Vert w_{0}-\beta
w\right\Vert $, $y_{2}\left(  \beta\right)  =\beta h-h_{0}$ and $y_{3}\left(
\beta\right)  =-r\left\Vert w\right\Vert \beta+r\operatorname{Re}\left\langle
w_{0},\tilde{w}\right\rangle .$}%
\end{figure}
\end{description}
\end{description}

\end{document}